\documentclass[DIV=22, paper=a4, fontsize=11pt]{article} 

\usepackage[english]{babel} 
\usepackage{amsmath,amsfonts,amsthm} 

\usepackage{fancyhdr} 
\usepackage{lastpage} 

\usepackage[bottom=1.5in,top=.7in,left=.9in,right=.9in]{geometry}

\fancyhfoffset[L]{0cm}
\fancyhfoffset[R]{0cm}

\usepackage{graphicx}      

\usepackage{caption}
\usepackage{subfigure}

\newcommand{\R}{\mathbb R}
\newcommand{\cN}{\mathcal N}

\newcommand{\cD}{\mathcal D}
\newcommand{\cY}{\mathcal Y}
\newcommand{\BQ}{\mathcal B_Q}
\newcommand{\cJ}{\mathcal J}

\newcommand{\C}{\mathcal C}
\newcommand{\cR}{\mathcal R}

\usepackage{titling} 

\newcommand{\HorRule}{ \rule{\linewidth}{1pt}} 

\pretitle{\vspace{-30pt} \begin{flushleft}  \bf\fontsize{22}{22} \HorRule\\ \selectfont} 

\title{A Successive Constraint Approach to Solving Parameter-Dependent Linear Matrix Inequalities} 

\posttitle{\par\end{flushleft}\vskip -0.5em} 

\preauthor{\begin{flushleft}\large \lineskip 0.5em \usefont{OT1}{phv}{b}{sl}} 

\author{Robert O'Connor, } 

\postauthor{\footnotesize \usefont{OT1}{phv}{m}{sl} 
  RWTH Aachen University 
  \par\end{flushleft}\HorRule}

\date{} 

\begin{document}

\maketitle 

\thispagestyle{fancy} 


\subsection*{Abstract:}
\textbf{We present a successive constraint approach that makes it possible to cheaply solve large-scale linear matrix inequalities for a large number of parameter values. The efficiency of our method is made possible by an offline/online decomposition of the workload. Expensive computations are performed beforehand, in the offline stage, so that the problem can be solved very cheaply in the online stage. We also extend the method to approximate solutions to semidefinite programming problems.}


\section{Introduction}
\label{}

Linear matrix inequalities (LMIs) are a general type of convex constraint that include linear as well as quadratic constraints \cite{VB2000,VB1996} and lead to very natural formulations of a large number of problems in control and systems theory \cite{BBF+1993,FAG1996}. They can be solved using a wide range of existing methods~\cite{AHO1998,FKM+2001,HR2000,Toh2003}, but that can be expensive for large-scale problems.  In particular, problems resulting from the discretization of partial differential equations can be extremely expensive due to their high dimensionality. The computations are even more expensive if parameter-dependent problems are considered and solutions are needed for a large number of parameter values.  In that case traditional solution methods are extremely inefficient.  We propose the construction of reduced-order models that take advantage of the parametric nature of the problem and allow us to very cheaply produce solutions for a large number of parameter values.

Let us introduce a finite-dimensional, bounded parameter domain $\cD\in\R^p$ and the parameter-dependent LMI
\begin{equation}\tag{L}\label{LMI}
F(x;\mu):=\sum_{q=1}^{Q_F}\left[\theta_q^0(\mu)+\theta_q^L(\mu)x\right]F_q\succeq 0.
\end{equation}
Here $F(x;\mu)\in\R^{\cN\times \cN}$ is a symmetric matrix that depends on both the parameter $\mu\in \cD$ and the decision variable $x\in\R^n$ and is composed of $Q_F$ parameter-independent matrices $F_q\in\R^{\cN\times \cN}$. The parameter dependencies of $F(x;\mu)$ are given by the functions $\theta^0(\cdot):\cD\rightarrow \R^{Q_F}$ and $\theta^L(\cdot):\cD\rightarrow \R^{Q_F\times n}$.  We use subscripts to indicate components of a vector, such that $\theta_q^0(\mu)$ is the $q$\textsuperscript{th} element of $\theta^0(\mu)$.  Similarly, we write $\theta_q^L(\mu)$ to indicate the $q$\textsuperscript{th} row of $\theta^L(\mu)$. The symbol $\succeq$ will be used in the sense that $P\succeq0$ indicates that the symmetric matrix $P$ is positive semi-definite.

The goal of this paper is to efficiently solve the following problems for a large number of parameter values $\mu\in\cD$:
\begin{enumerate}
\item The strict feasibility problem: Find an $x\in\R^n$ such that $F(x;\mu)\succ0$.
\item The semidefinite program (SDP):
\begin{equation}\tag{S}\label{SDP}
 \underset{x\in\R^n}{\text{minimize}}\quad
 c(\mu)^Tx \quad
 \text{subject to}\quad
  F(x;\mu)\succeq 0.
\end{equation}
\end{enumerate}
Rather than directly solving these problems for each new parameter value, we will solve them for only a small set of intelligently chosen parameter values. We will then use the resulting solutions to build a reduced-order model that can approximate the solution anywhere in $\cD$.

The method that we propose can be viewed as a generalization of the successive constraint method (SCM) \cite{CHM+2008,HRS+2007}, which is often used in the field of reduced basis methods to evaluate stability constants \cite{RHP2008}. This method will allow us to very cheaply determine feasible solutions for any $\mu\in\cD$.  That is made possible by decomposing the computational workload into offline and online stages.  All expensive computations will be performed in advance, during the offline stage.  During the online stage the cost to solve the problem for a new parameter value will be independent of the size of the original constraint, $\cN$. In that way the computational cost of each new solution will remain cheap even if the original constraint has very large dimensions. 

Our methods are applicable to a wide range of LMIs and can also be used to extend the applicability of SCM. In the context of reduced basis methods, applications could involve bounding stability constants with respect to parameter-dependent norms or the selection of Lyapunov functions for the computation of error bounds \cite{OConnor2016}.  In Section \ref{sec:example} we present an example in which we optimize a system while ensuring that it remains stable.


\section{Reduced-order modeling for strict feasibility}\label{sec:feas}

SCM was originally designed to approximate coercivity constants. We will apply a modified version of SCM to the coercivity constant
\begin{equation}
\alpha(x;\mu):=\inf_{v\in \R^\cN}\frac{v^TF(x;\mu)v}{v^TF_Sv},
\end{equation}
where $F_S\in\R^{\cN\times \cN}$ is a fixed symmetric positive-definite matrix. From the definition it is clear that $\alpha(x;\mu)\geq0$ is equivalent to $F(x;\mu)\succeq0$ for all symmetric positive definite matrices $F_S$.  Nevertheless, an appropriate choice of $F_S$ could be beneficial from a numerical point of view. If we are dealing with PDE discretizations, it can be advantageous to choose a matrix associated with an energy norm.

The first step in applying SCM is reformulating the coercivity constant as follows:
\begin{equation}\label{eq:coer_opt}
\alpha(x;\mu)=\inf_{y\in \cY}\left[\theta^0(\mu)+\theta^L(\mu)x\right]^Ty,\quad\text{where}\quad
\cY:=\left\{y\in\R^{Q_F}\Big|y_q=\frac{v^TF_qv}{v^TF_Sv},v\in\R^\cN\right\}.
\end{equation}
This formulation has the advantage that the complexity of the problem has been shifted to the definition of the set $\cY$.  That allows us to compute lower and upper bounds for $\alpha(x;\mu)$ by approximating $\cY$.

A lower bound for $\alpha(x;\mu)$ can be derived by approximating $\cY$ from the outside. A bounded but primitive approximation for $\cY$ is given by
\begin{equation}
\BQ:=\prod_{i=1}^{Q_F} \left[\inf_{y\in\cY}y_q,\sup_{y\in\cY}y_q\right]\subset \R^{Q_F}.
\end{equation}
$\BQ$ will generally be much larger than $\cY$ so we will add constraints to restrict it and better approximate $\cY$. Let us assume that we have computed lower bounds $\bar\alpha(\bar x;\bar\mu)\leq\alpha(\bar x;\bar \mu)$ for each pair $(\bar x,\bar \mu)$ in a predetermined set $\cR\subset\R^n\times\cD$. The set $\cY$ is then contained in the set
\begin{equation}
\cY_{out}:=\left\{y\in\BQ\Big|\left[\theta^0(\bar\mu)+\theta^L(\bar\mu)\bar x\right]^Ty\geq\bar\alpha(\bar x;\bar\mu),\, \forall\ (\bar x,\bar\mu)\in \cR\right\}.
\end{equation}
Using $\cY_{out}$ we can define
\begin{equation}\label{eq:coer_lb}
\alpha_{out}(x;\mu):=\inf_{y\in \cY_{out}}\left[\theta^0(\mu)+\theta^L(\mu)x\right]^Ty,
\end{equation}
such that $\alpha_{out}(x;\mu)\leq\alpha(x;\mu)$ for any $x\in\R^n$ and any $\mu\in\cD$.

The lower bound $\alpha_{out}(x;\mu)$ can be formulated as the solution to a linear programming problem.  To see that we note that $\cY_{out}$ is a linearly constrained subset of $\R^{Q_F}$. As a result we can find a matrix $A_{out}\in\R^{\ell\times Q_F}$ and a vector $b_{out}\in\R^\ell$ such that $\cY_{out}=\{y\in\R^{Q_F}|A_{out}y\geq b_{out}\}$. Here the symbol $\geq$ indicates a componentwise comparison of vectors and $\ell\ll\cN$ is the number of constraints.  The lower bound $\alpha_{out}(x;\mu)$ can be calculated, for any $x\in\R^n$ and any $\mu\in\cD$, as the solution to the linear program (\ref{P}) or its dual (\ref{D}):

\noindent
\fbox{\begin{minipage}{0.482\textwidth}
\vspace{.1cm}
\begin{equation}\tag{P}\label{P}
\begin{aligned}
& \underset{y\in\R^{Q_F}}{\text{minimize}}
& & \left[\theta^0(\mu)+\theta^L(\mu)x\right]^Ty \\
& \text{subject to}
& &  A_{out}y\geq b_{out}
\end{aligned}
\end{equation}
\vspace{.1cm}
\end{minipage}}
\fbox{\begin{minipage}{0.482\textwidth}
\begin{equation}\tag{D}\label{D}
\begin{aligned}
& \underset{p\in\R^\ell}{\text{maximize}}
& & b_{out}^Tp \\
& \text{subject to}
& &  A_{out}^Tp=\theta^0(\mu)+\theta^L(\mu)x\\
& & & p\geq0
\end{aligned}
\end{equation}
\end{minipage}}

To find an optimal value of $x$ we will maximize (\ref{D}) over all $x\in\R^n$ using the problem
\begin{equation}\tag{RF}\label{RF}
\underset{p\in\R^\ell,x\in\R^n}{\text{maximize}}\quad b_{out}^Tp\quad \text{subject to}
\quad  A_{out}^Tp-\theta^L(\mu)x=\theta^0(\mu) \quad \text{and} \quad p\geq0.
\end{equation}
If $b_{out}^Tp$ is strictly greater than $0$, then the associated $x$ is guaranteed to strictly satisfy (\ref{LMI}). Even if that is not the case, (\ref{RF}) is always feasible. To show that we note that the boundedness of $\cY_{out}\subset\BQ$ implies the boundedness of (\ref{P}). By duality (\ref{D}) is then feasible for any $x\in\R^n$ and any $\mu\in\cD$ which implies that (\ref{RF}) is feasible for any $\mu\in\cD$.  (\ref{RF}) is also bounded if $\max_{x\in\R^n} \alpha(x;\mu)$ is bounded; otherwise, we will be content with any $x$ associated with a large, positive value of $b_{out}^Tp$. The main reason to use (\ref{RF}) is that it is cheap to solve: It is a small linear program with $n+\ell$ variables, $Q_F$ equality constraints, and $\ell$ inequality constraints.

Since it is not possible to work directly with the infinite set $\cD$, we introduce a large set $\Xi\subset\cD$ of discrete points that are representative of $\cD$ and a much smaller set $\C_k\subset\Xi$ with cardinality $k$. The set $\cR$ will be made up of parameter values $\bar\mu\in\Xi$ and associated approximate solutions $\bar x\in\R^n$.  For parameter values $\bar\mu$ that are in $\C_k$ we will use solutions $\bar x$ to the original problem (\ref{LMI}) and set $\bar\alpha(\bar x;\bar\mu)=\alpha(\bar x;\bar\mu)$.  For parameter values in the much larger set $\Xi\setminus\C_k$ we will use solutions $\bar x$ to (\ref{RF}), and set $\bar\alpha(\bar x;\bar\mu)=\alpha_{out}(\bar x;\bar\mu)$. The pairs $(\bar x,\bar\mu)$ that will be included in $\cR$ will depend on the parameter $\mu$ for which we would like a feasible solution. For two natural numbers $M_C$ and $M_\Xi$ we define $\cR$ to be a set of $M_C$ pairs $(\bar x,\bar\mu)$ with $\bar\mu\in\C_k$ and $M_\Xi$ pairs $(\bar x,\bar\mu)$ with $\bar\mu\in\Xi\setminus\C_k$.  In both cases we will choose the pairs $(\bar x,\bar\mu)$ with $\bar\mu$ closest to $\mu$ in a predetermined metric. If $M_C>k$, we simply include all of $\C_k$ in $\cR$.


To build our model we will make use of the greedy algorithm that was introduced in the context of reduced basis methods~\cite{VPR+2003} and plays a vital role in SCM. The algorithm is initiated by choosing a small initial set $\C_k$. For each $\bar\mu\in\C_k$ the problem (\ref{LMI}) is solved once to determine the associated values of $\bar x$ and $\bar\alpha(\bar x;\bar\mu)$. For $\bar\mu\in\Xi\setminus \C_k$ we initially set $\bar x=0$ and $\bar \alpha(\bar x;\bar\mu)=-\infty$. The model is then improved iteratively. In each iteration we will solve (\ref{RF}) for each $\bar\mu\in\Xi\setminus\C_k$ and update the stored values of $\bar x$ and $\bar \alpha(\bar x;\bar \mu)$.  The $\bar\mu$ associated with the smallest value of $\bar\alpha(\bar x;\bar\mu)$ is then added to $\C_k$, and (\ref{LMI}) is solved to update the stored values of $\bar x$ and $\bar\alpha(\bar x;\bar\mu)$. The process terminates when the smallest value of $\bar\alpha(\bar x;\bar\mu)$ is larger than some positive tolerance.

Solving problems with our method involves two stages.  During the offline stage the greedy algorithm builds the model.  The expensive operations in the offline stage are $k$ solutions of (\ref{LMI}) and nearly $k|\Xi|$ solutions of (\ref{RF}). Here $k$ is the cardinality of the final $\C_k$ and $|\Xi|$ is the cardinality of $\Xi$.  During the online stage (\ref{RF}) can be constructed and solved very cheaply for any new parameter value. The online computational cost is independent of $\cN$ to ensure that it remains cheap even if $\cN$ is very large.

During the greedy algorithm it is necessary to solve (\ref{LMI}) for various parameter values.  That can be done using a variety of methods including algorithms that solve semidefinite programming problems \cite{AHO1998,FKM+2001,HR2000,Toh2003}.  In this section we propose a method that is based on SCM and allows us to reuse information that we have already computed.

For a predetermined, finite subset $\cY_{in}$ of $\cY$ we define
\begin{equation}\label{eq:coer_ub}
\alpha_{in}(x;\mu):=\inf_{\bar y\in \cY_{in}}\left[\theta^0(\mu)+\theta^L(\mu)x\right]^T\bar y,
\end{equation}
such that $\alpha_{in}(x;\mu)\geq\alpha(x;\mu)$ for all $x\in\R^n$ and $\mu\in\cD$. Here we will make use of a variation of the greedy algorithm.  In each iteration we solve $\max_{x\in\R^n}\alpha_{in}(x;\mu)$ to find the optimal $x$. For that value of $x$ we then compute $\alpha(x;\mu)$ and an element $y$ of $\cY$ for which the infimum in (\ref{eq:coer_opt}) is reached.  That $y$ is then added to $\cY_{in}$ to improve the upper bound $\alpha_{in}(\cdot;\mu)$. The algorithm terminates when $\alpha(x;\mu)$ is sufficiently large or sufficiently close to $\alpha_{in}(x;\mu)$.

In each iteration of the algorithm both the small linear program associated with $\max_{x\in\R^n}\alpha_{in}(x;\mu)$ and the $\cN$-dimensional eigenvalue problem associated with $\alpha(x;\mu)$ need to be solved once.  The dominant cost is that of the eigenvalue problems.  To reduce the number of iterations and eigenvalue solves we will store the last value of $\bar y$ that is calculated each time we solve (\ref{LMI}) for a new $\mu\in\C_k$.  The set $\cY_{in}$ can then be initialized using those values of $\bar y$.

\section{The semidefinite programming problem}\label{sec:SDP}

In this section we will show how our methods can be used to approximate solutions to the semidefinite program (\ref{SDP}). The idea is to minimize $c(\mu)^Tx$ over all $x\in\R^n$ that satisfy $\alpha_{out}(x;\mu)\geq0$.  That is done using the following problem:
\begin{equation}\tag{RS}\label{RS}
\underset{x\in\R^n,p\in\R^\ell}{\text{minimize}}\quad c(\mu)^Tx \quad\text{subject to}\quad
\begin{bmatrix}A_{out}^T& -\theta^L(\mu)\end{bmatrix}
\begin{bmatrix}p\\x\end{bmatrix}=\theta^0(\mu)\quad \text{and}\quad  \begin{bmatrix}I\\b_{out}^T\end{bmatrix}p\geq
\begin{bmatrix}0 \\ 0\end{bmatrix}.
\end{equation}
Here $A_{out}$, $b_{out}^T$ and $\alpha_{out}(x;\mu)$ are all the same as in Section \ref{sec:feas} except that the stored values of $\bar x$ will be calculated using (\ref{SDP}) and (\ref{RS}) rather than (\ref{LMI}) and (\ref{RF}). The new constraint $b_{out}^Tp\geq0$ guarantees that $\alpha_{out}(x;\mu)\geq0$.

Let us write $\cJ_{out}$ and $\cJ$ to denote the optimal values of (\ref{RS}) and (\ref{SDP}), respectively. Based on the fact that the optimal $x\in\R^n$ for (\ref{RS}) will be feasible for (\ref{SDP}) we know that $\cJ_{out}\geq \cJ$. Ideally $\cJ_{out}$ will also be a good approximation to $\cJ$ despite being much cheaper to compute: It requires the solution of a linear program with $n+\ell$ variables, $Q_F$ equality constraints, and $\ell+1$ inequality constraints.

For any $\mu\in\cD$, the boundedness of (\ref{SDP}) implies the boundedness of (\ref{RS}), but the question of feasibility is more complicated.  To ensure feasibility we first build a model (\ref{RF}) to find strictly feasible solutions.  We then convert that model to the form of (\ref{RS}).  In doing so we consider the stored values of $\bar x\in\R^n$, which should be feasible, to be approximate solutions to (\ref{SDP}) and we reinitialize $\C_k$ as an empty set.  We can then improve the model using another greedy algorithm.

The goal of this greedy algorithm is to reduce the error $\cJ_{out}-\cJ$. Since it is too expensive to compute $\cJ$ online, we will compute a lower bound for it by considering the problem
\begin{equation}\tag{ER}\label{ER}
 \underset{x\in\R^n}{\text{minimize}}\quad
 c(\mu)^Tx \quad
 \text{subject to}\quad
  \alpha_{in}(x;\mu)\geq0.
\end{equation}
Here $\alpha_{in}(x;\mu)$ is defined like in (\ref{eq:coer_ub}) with $\cY_{in}$ initially being the set of $\bar y\in\cY$ associated with parameter values in $\C_k$. The optimal value of (\ref{ER}), which we will denote $\cJ_{in}$, is by definition a lower bound for $\cJ$.  That allows us to bound the error $\cJ_{out}-\cJ\geq0$ from above using $\cJ_{out}-\cJ_{in}$.  Like the previous approximations, this bound can be computed online with a number of computations that is independent of $\cN$. Although (\ref{ER}) is always feasible, it could be unbounded.  That can be fixed by increasing $\C_k$ and consequently $\cY_{in}$.

During each iteration of the greedy algorithm we will solve both (\ref{RS}) and (\ref{ER}), compute the difference $\cJ_{out}-\cJ_{in}$, and update both the estimates $\bar x\in\R^n$ and the associated values $\bar\alpha(\bar x;\bar\mu)$ for each $\bar\mu\in\Xi\setminus\C_k$.  The $\bar\mu$ that produced the maximum value of $\cJ_{out}-\cJ_{in}$ is then added to $\C_k$ and a more accurate solution to (\ref{SDP}) is computed. That allows us to update both (\ref{RS}) and (\ref{ER}).  This process terminates when the maximum value of $\cJ_{out}-\cJ_{in}$ is below a desired tolerance.

The major computational burdens in the offline stage of our method are the construction of the model (\ref{RF}), $k$ solutions of (\ref{SDP}), and nearly $k|\Xi|$ solutions of both (\ref{RS}) and (\ref{ER}).  During the online stage an approximate solution can be determined for any given $\mu\in\cD$ by simply constructing and solving (\ref{RS}).  The error can also be bounded online by solving (\ref{ER}).

The more accurate solutions that we need can be computed using either general solvers for semidefinite programming problems \cite{AHO1998,FKM+2001,HR2000,Toh2003} or a method that is similar to the one presented in Section \ref{sec:feas}.  For the latter we will consider the following problem with a fixed value of $\mu\in\cD$:
\begin{equation}\tag{T}\label{T}
 \underset{x\in\R^n}{\text{minimize}}\quad
 c(\mu)^Tx \quad
 \text{subject to}\quad
  \alpha_{in}(x;\mu)>\alpha_{min},
\end{equation}
where $\alpha_{min}>0$ is some small constant that is used as a tolerance. We again use a sort of greedy algorithm. In each iteration we solve (\ref{T}) to get the approximate solution $x$ and update $\cY_{in}$ by adding the new value of $y$. The process is repeated until the stopping condition $\alpha(x;\mu)>0$ is satisfied. If $\alpha_{min}$ is sufficiently small, the resulting $x$ should be a good approximation to the true solution.

\section{Numerical example: System stabilization}\label{sec:example}

\begin{figure}
 \centering
 \subfigure[The domain of the problem]
 {
  \includegraphics[width=3cm,trim=0cm 0cm 0cm 0cm]{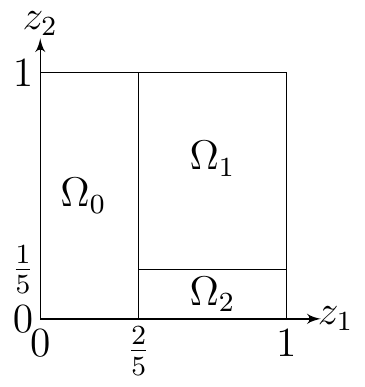}
   \label{desc:dom}
   }
 \subfigure[Feedback loop: $u(t)=-x\zeta(t)$]
 {
  \includegraphics[width=6.5cm,trim=0cm 0cm 1cm 0cm]{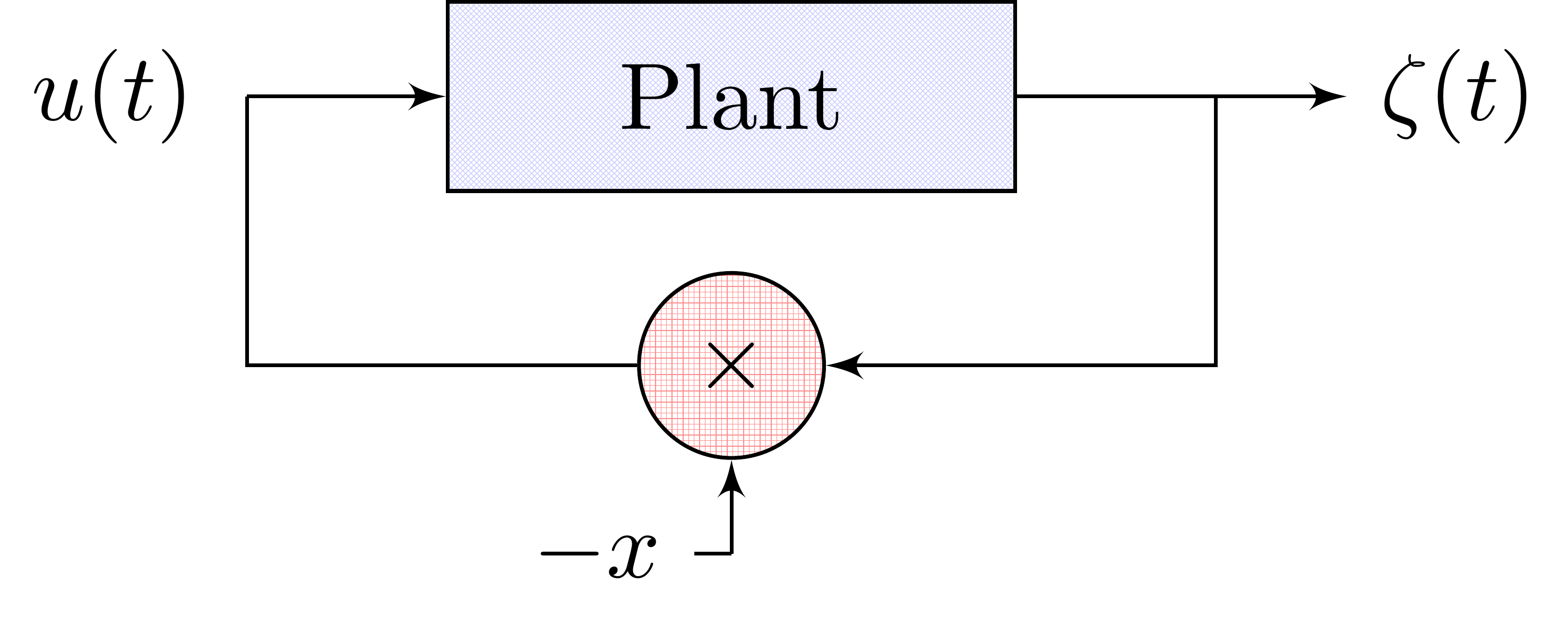}
   \label{desc:fb}
   }

 \label{desc}
 \caption{Setup of the numerical example}\vspace{-.2cm}
\caption*{Configuration de l'exemple num\'erique}
\vspace{-.1cm}
\end{figure}

To test our method we will consider a reaction-diffusion equation on the unit square, $\Omega$, depicted in Figure \ref{desc:dom}.  To facilitate the explanation we consider a problem with just one parameter $\mu\in\cD:=[0,3]$ and one decision variable $x\in\R$.  The problem is related to the following system:
\begin{equation}
\dot y=\Delta y+\mu y \mathbf 1_{\Omega_1}+\mathbf 1_{\Omega_2}u,\quad \forall z\in\Omega;\qquad\frac{\partial y}{\partial \eta}=0, \quad \forall z\in\partial \Omega;\qquad u=-x\zeta;\qquad \zeta=\int_{\Omega_2}y.
\end{equation}

\noindent Here $\mathbf 1_{S}$ indicates the characteristic function of a set $S$ and $\eta$ is the outward-pointing unit vector on the boundary $\partial \Omega$. The decision variable $x$ determines the feedback gain from the system's output $\zeta(t)$ to its input $u(t)$ as depicted in Figure \ref{desc:fb}.  We use linear finite elements to discretize the problem in space.  That gives us the following $2601$-dimensional closed-loop semi-discrete system
\begin{equation}\label{eq:semi_disc}
M\dot y+A(x;\mu)y=0, \quad\text{where}\quad A(x;\mu)=A_0-\mu A_1+xA_2,\quad \text{with}\quad M,A_0,A_1,A_2\in\R^{2601\times 2601}.
\end{equation}
Here $M$ is symmetric positive definite and $A(x;\mu)$ is symmetric.  Under these assumptions the system (\ref{eq:semi_disc}) is strictly stable for a given $\mu\in\cD$ and $x\in\R$ iff $A(x;\mu)\succ0$.  Our goal will be to minimize the cost of the control, by keeping $x$ small, while ensuring strict stability. Noting that $A_0$, which is associated with the operator $-\Delta$, is positive semidefinite and that $A_0+A_1$ is positive definite, we define $F(x;\mu):=(1-\rho)A_0+(-\mu-\rho)A_1+xA_2$ and $F_S:=A_0+A_1\succ0$ for some small $\rho>0$.  It then holds that
\begin{equation}
F(x;\mu)\succeq0\qquad\Longleftrightarrow \qquad A(x;\mu)- \rho F_S\succeq0,
\end{equation}
and hence the system in (\ref{eq:semi_disc}) is strictly stable if $F(x;\mu)\succeq0$.  We will search for a minimal stabilizing gain $x\in\R$ as the solution to (\ref{SDP}).

For our model we will use values of $M_C=4$ and $M_\Xi=3$. $\C_k$ is initialized using the smallest and the largest values in $\cD$. For this particular problem that guarantees that (\ref{RS}) is feasible for all $\mu\in\cD$. As a result we can start directly with (\ref{RS}) and do not need to build the model given in (\ref{RF}). Another modification of our method is that we will work with the relative error given by $(\cJ_{out}-\cJ_{in})/\cJ_{in}$ rather than the error $\cJ_{out}-\cJ_{in}$. That is possible because $\cJ_{in}$ will always be strictly positive.

\begin{figure}
 \centering
 \subfigure[Convergence of the full-order solver]
 {
  \includegraphics[width=4.7cm,trim=1cm 1.3cm 1cm 0cm]{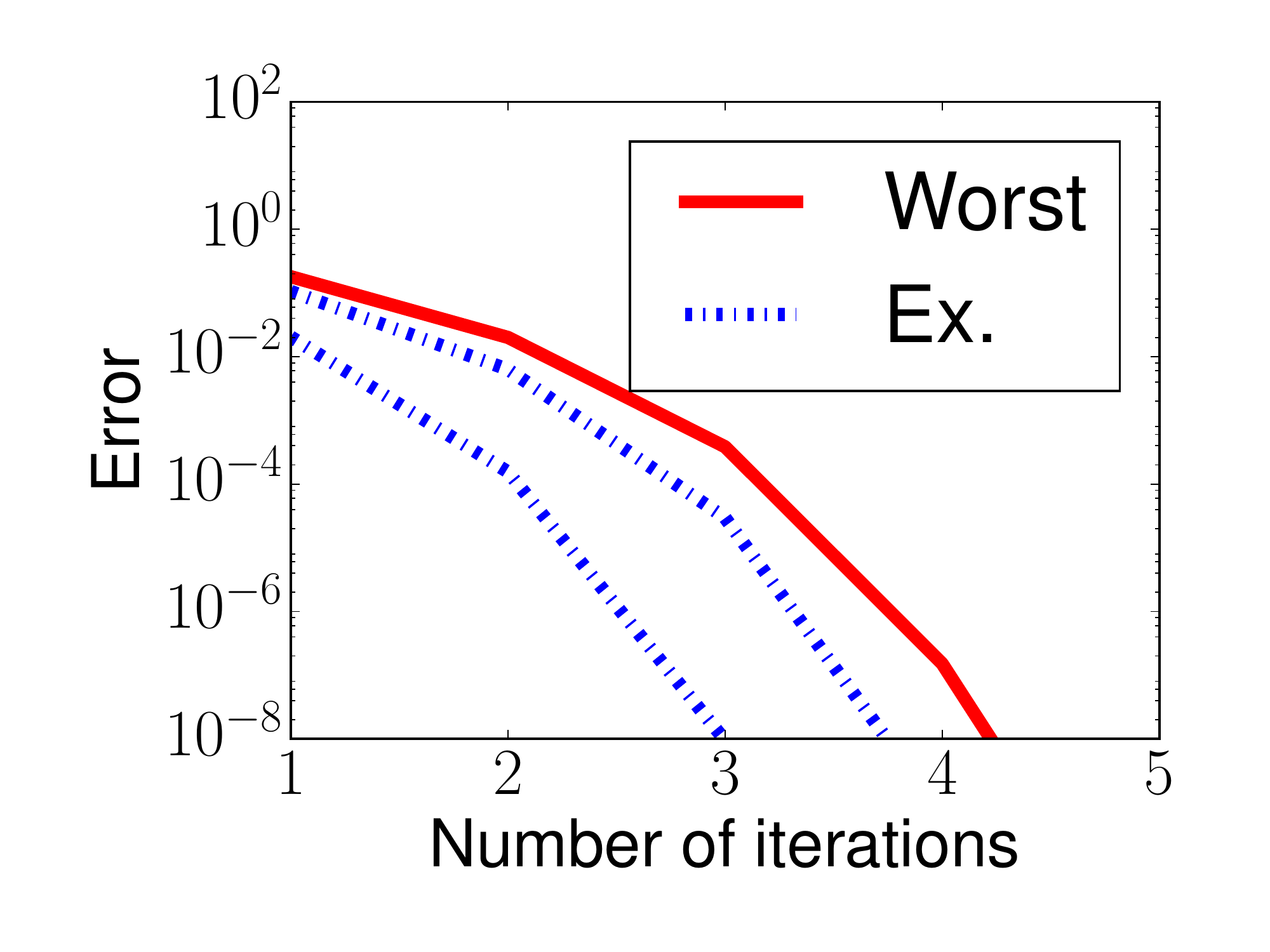}
   \label{res:fixed_mu}
   }
 \subfigure[Greedy Convergence]
 {
  \includegraphics[width=4.7cm,trim=1cm 1.3cm 1cm 0cm]{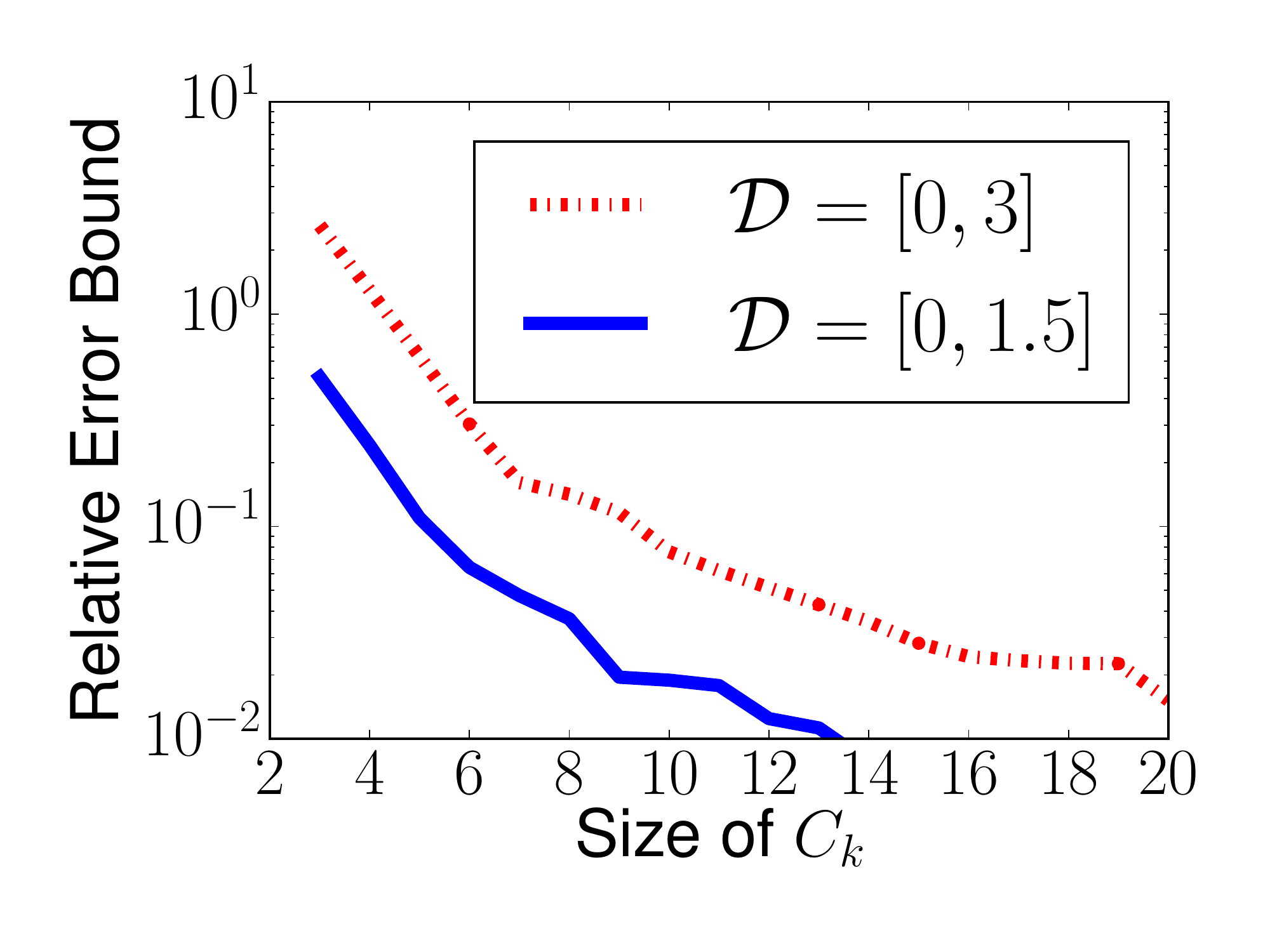}
   \label{res:greedy}
   }
   \subfigure[Uniform distribution vs. greedy]
 {
  \includegraphics[width=4.7cm,trim=1cm 1.3cm 1cm 0cm]{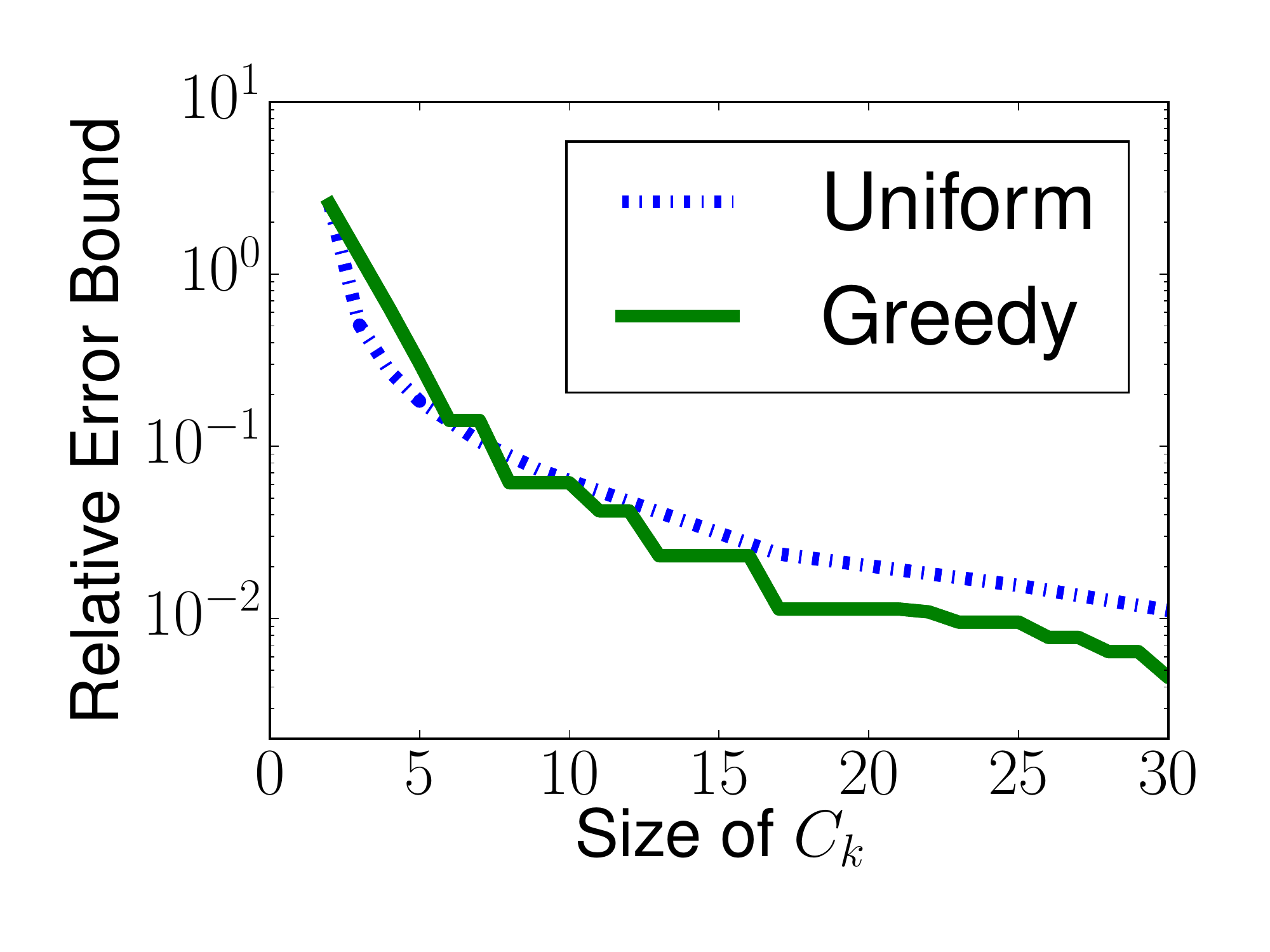}
   \label{res:snaps}
   }
\vspace{-.3cm}

 \label{res}
 \caption{Results of the numerical experiments}\vspace{-.2cm}
\caption*{R\'esultats des esp\'eriences num\'eriques}
\vspace{-.1cm}
\end{figure}

Figure \ref{res:fixed_mu} shows the convergence of the solver described in Section \ref{sec:SDP} for the computation of more accurate solutions to (\ref{SDP}). Here the error is measured as $\alpha_{in}(\bar x;\bar\mu)-\alpha(\bar x;\bar\mu).$  Two examples are shown as well as the worst case over $30$ random parameter values.  We recall that the most expensive part of this method is the eigenvalue solves and that one such solve is needed for each iteration. Figure \ref{res:greedy} shows the convergence of the greedy algorithm for two different parameter domains $\cD$.  The plotted values are the worst relative errors over the respective sets $\Xi$. For $\cD=[0,3]$ we set $\Xi$ to be $300$ uniformly distributed points and for $\cD=[0,1.5]$ we use $150$ points. Figure \ref{res:snaps} compares the results when $\C_k$ is constructed using the greedy algorithm or simply as uniformly distributed points. The plot shows the worst relative errors over $100$ random values of $\mu\in\cD$. The greedy algorithm performs slightly better than uniform distributions when $\C_k$ is large and also has the advantage that it is iterative.

{\small
  \setlength{\parskip}{0pt}
  \setlength{\itemsep}{0pt plus 0.1ex}
  \bibliographystyle{abbrv}
  \bibliography{../../../References/references.bib}
}

\end{document}